\documentclass[12pt]{article}
\usepackage{amsfonts,amssymb,amsthm}

\def\be{\begin{equation}}
\def\ee{\end{equation}}
\def\bea{\begin{eqnarray}}
\def\eea{\end{eqnarray}}

\def\m{\mu}

\def\tfrac#1#2{{\scriptstyle{\frac{#1}{#2}}}}         

\begin{document}
\thispagestyle{empty}
\hfill\today
\vspace{2.5cm}

\begin{center}
{\LARGE{\bf{Quantization of 
Drinfel'd doubles}}}
\end{center}

\bigskip\bigskip

\begin{center}
A. Ballesteros$^1$, E. Celeghini$^2$  and M.A. del Olmo$^3$
\end{center}

\begin{center}
$^1${\sl Departamento de F\'{\i}sica, Universidad de Burgos, \\
E-09006, Burgos, Spain.}\\
\medskip

$^2${\sl Departimento di Fisica, Universit\'a  di Firenze and
INFN--Sezione di
Firenze \\
I50019 Sesto Fiorentino,  Firenze, Italy}\\
\medskip

$^3${\sl Departamento de F\'{\i}sica Te\'orica, Universidad de
Valladolid, \\
E-47011, Valladolid, Spain.}\\
\medskip

{e-mail: angelb@ubu.es, celeghini@fi.infn.it, olmo@fta.uva.es}
\end{center}

\bigskip

\begin{abstract}
Hopf algebra quantizations of 4-dimensional and 6-dimensional real classical
Drinfel'd doubles are studied by following a direct ``analytic" approach. The full
quantization is explicitly obtained for most of the Drinfel'd
doubles, except a small number of them for which the dual Lie algebra is either
$\mathfrak{sl}_2$ or
$\mathfrak{so}(3)$. In the latter cases, the classical $r$-matrices underlying the 
Drinfel'd double quantizations contain  known standard ones plus additional twists.
Several new four and six-dimensional quantum algebras are presented and some general
features of the method are emphasized.
\end{abstract}
\vskip 1cm

\vskip 0.4cm

\vfill
\eject
\section{Introduction\label{introduccion}}

The essential role that Lie bialgebras
play in the quantization of Poisson-Lie structures
comes from the well-known result by Drinfel'd \cite{Dr} that
establishes a one to one correspondence between Poisson-Lie groups
and Lie bialgebras. In fact, the concept of
classical double \cite{Dr} is just a reformulation of that of  Lie  bialgebra
in terms of a (double) dimensional Lie algebra  endowed with a suitable pairing. 
This ``duplication"
process can be iterated by taking into account that the double Lie algebra  can be in turn
equipped with a quasitriangular Lie bialgebra structure  by means of a canonical classical
$r$-matrix.

The Hopf algebra quantization of such a double
Lie bialgebra is the so-called quantum double, a basic object in quantum group theory 
(see, for instance, \cite{Dri}--\cite{contr} for a detailed exposition and references
therein). In particular, quantum doubles are helpful for the explicit construction of
quantum
$R$-matrices of quantum groups and supergroups by following a ``universal ${\cal
T}$-matrix approach" \cite{FG} (see 
\cite{Ge}-\cite{BRST}) and have been also considered as symmetries in quantum field
theory \cite{QFT,NPB}. Moreover, $\sigma$-models related by Poisson-Lie T-duality are
directly connected with classical doubles (see \cite{duality}-\cite{Snobl} and references
therein). 

We also recall that in Ref.~\cite{gomez} the full classification
of two and three dimensional (3D) real Lie bialgebra structures and their associated
Drinfel'd doubles has been obtained. That result is tantamount to the (first order)
classification of quantum deformations of 2D and 3D real Lie algebras. Complementarily,
the classification of 4D and 6D Drinfel'd
doubles has been performed in \cite{Snobl} and \cite{HS},  respectively, by following a
``direct" approach. Throughout the paper we will preserve the notation and labeling for
the Lie bialgebras and Drinfel'd doubles given in
\cite{gomez}, although both classifications are fully equivalent (up to the fact that
\cite{gomez} does not consider 6D classical doubles coming from a trivial Lie bialgebra
with vanishing cocommutator $\delta$ or, equivalently, with abelian dual
$\mathfrak{g}^\ast$).

The aim of this work is to provide a global insight into
the Hopf algebra quantization of 4D and 6D Drinfel'd doubles (DD) by making use
of a direct quantization procedure that has been recently used in order to obtain and
classify 3D quantum algebras \cite{3dim}.  The motivation for this study is
two-fold. 

Firstly, the family of 6D DD algebras is physically interesting. For
instance, it contains the $so(3,1)$ and $so(2,2)$
algebras as well as the (2+1) Poincar\'e algebra. Therefore the DD
quantizations will provide quantum deformations for these
algebras. In particular, the class of  DDs with $\mathfrak{g}=\mathfrak{r}_3(1)$ can be
thought as $so(p,q)$ Lie algebras (with $p+q=4$) and some of their contractions
\cite{so22}.  In general, such quasitriangular quantizations turn out to be
superpositions of a standard (quasitriangular) quantization plus twists. We will discuss
the properties of all these quantum algebras and we explicitly obtain
some of them by making use of a  quantization
procedure recently introduced by us \cite{3dim}.

Secondly, since Lie bialgebras (therefore, DD algebras)
provide a classification of first-order quantum deformations, the set of DD
quantizations can be considered as an appropriate setting for the
search of classification schemes for quantum deformations. In fact, several
strong regularities are found within the set of quantum algebras that are presented in
this work, and these common facts can be thought of guidelines for a future research
program of a classification ``\'a la Cartan" of quantum groups \cite{3dim}. Among them,
we mention the following common properties of the quantizations:

\begin{itemize}

\item All 4D and 6D DD algebras are non-simple Lie algebras.

\item The only functions appearing in the deformed commutation rules
and coproducts are exponentials and polynomials \cite{Varadarajan}.

\item The quantization can be obtained by following a ``analytic
approach" in the full
symmetrized basis of the quantum universal enveloping algebra that preserves a
``generalized cocommutativity" property \cite{3dim}. 

\item All the DD deformations coming from the canonical skew-symmetric DD $r$-matrix
are standard ones. 

\item The ``crossed" commutation rules between the algebra and the dual contain the most
complex part of the quantum deformations. 

\end{itemize}

The structure of the paper is as follows: Section 2 is devoted to fix the notation.
Section 3 presents the study of  quantum 2D DD algebras as a toy-model for
further quantizations.  The following sections develop the quantizations of
6D DD algebras by following the classification given in \cite{gomez}. Finally,
some remarks close the paper.

\section{Drinfel'd Double (bi)algebras}

Let us consider
a Lie bialgebra ($\mathfrak{g},\delta$) and
a basis $\{x^i\}$  of $\mathfrak{g}$. Such a Lie
bialgebra can be characterized
 by a pair of structure tensors
($f^{lm}_n,c^k_{ij}$), i.e., 
$$ [x^i,x^j]=f^{ij}_k x^k, \qquad \delta(x^n)=c^{n}_{lm} x^l\otimes x^m.  
\label{agc}
$$ In this language, the cocycle condition for the cocommutator $\delta$
becomes the following compatibility condition between the
tensors $c$ and $f$
\be
f^{ab}_k c^k_{ij} = f^{ak}_i c^b_{kj}+f^{kb}_i c^a_{kj}
+f^{ak}_j c^b_{ik} +f^{kb}_j c^a_{ik}. \label{agb}
\ee  

Now we fix a basis $\{X_i\}$ for the dual algebra $\mathfrak{g} ^*$ through the
following pairing 
$$  \langle X_i,X_j\rangle= 0,\quad \langle x^i,x^j\rangle=0, \quad
\langle x^i,X_j\rangle= \delta^i_j,\quad \forall i,j.\label{age}
$$ Then ($\mathfrak{g} ^*,\eta$) is also a Lie
bialgebra with structure tensors ($f, c$), i.e.,
$$ [X_i,X_j]=c^k_{ij}X_k, \qquad \eta(X_n)= f^{lm}_n X_l\otimes X_m.  
\label{aga}$$ 

This duality leads to the consideration of
the pair ($\mathfrak{g} ,\mathfrak{g}^*$) and its associated vector space 
$\mathfrak{g} \oplus \mathfrak{g} ^*$, that can be endowed with a Lie algebra
structure by means of the commutators
\be
[x^i,x^j]= f^{ij}_k x^k, \quad
[X_i,X_j]= c^k_{ij}X_k, \quad  
[x^i,X_j]= c^i_{jk}x^k- f^{ik}_j X_k.\label{agd}
\ee
This Lie algebra, $D ({\mathfrak{g}})$, is called the
Double Lie algebra of $(\mathfrak{g},\delta)$. Obviously,
$\mathfrak{g}$ and $\mathfrak{g}^*$ are subalgebras of $D ({\mathfrak{g}})$, and the 
compatibility condition (\ref{agb}) is just the Jacobi identity
for (\ref{agd}).

Moreover, if $\mathfrak{g}$ is a finite dimensional Lie algebra,
then $D ({\mathfrak{g}})$  can be endowed with a
(quasitriangular) Lie bialgebra structure $(D ({\mathfrak{g}}),\delta_{DD})$  generated
by the classical
$r$-matrix

\be
r=\sum_i{x^i\otimes X_i}
\label{rmat}
\ee
or, equivalently, by its skew-symmetric counterpart 
$$ \tilde r=\frac12 \sum_i{x^i\wedge X_i}.
\label{rmatskew}
$$ In this respect, note that
$$ C=\sum_i{(x^i\,X_i+X_i\,x^i)} 
\label{cas}
$$ is always a Casimir operator for the DD algebra. Hence, if we denote
$$ \Omega=\sum_i{(x^i\otimes X_i+X_i\otimes x^i)},
$$ then
$[1\otimes Y + Y
\otimes 1,\Omega]=0$ for any generator $Y$ of the DD
algebra and $
\tilde r=r - \frac12 \Omega
$.

The cocommutator $\delta_{DD}$  derived from (\ref{rmat}) is
$$ \delta_{DD}(x^i)=\delta(x^i)=c^i_{jk}\,x^j\otimes x^k\qquad
\delta_{DD}(X_i)=-\eta(X_i)=- f_i^{jk}\,X_j\otimes
X_k.
\label{codob}
$$ In fact this ``double Lie bialgebra" has as sub-Lie-bialgebras the 
original one $(\mathfrak{g},\delta)$ and its dual ($\mathfrak{g}^\ast
,\eta$).

\section{Quantum two-dimensional DD algebras}

In this section we shall illustrate the procedure that we are going to follow
in the rest of the paper by considering the quantization of the 2D DD
 algebras. In general, we shall proceed as follows:
\begin{enumerate}
 \item
 We obtain a quantum coproduct and deformed commutation rules for each of the
Hopf subalgebras that quantize the sub-Lie-bialgebras $(\mathfrak{g},\delta)$ and 
($\mathfrak{g}^\ast ,\eta$). 

\item
 We find the deformed counterparts for the ``crossed" commutation rules
$[x^i,X_j]$.
\end{enumerate}

We will use the Gomez results and notation for Lie bialgebras \cite{gomez},
that in the two-dimensional case are shown in Table I. In the first row the Lie
sub-Lie-bialgebra  structure of 2D DD algebras is described, and their Lie
algebra brackets are explicitly given in the most appropriate basis. 

\vfill\eject

\centerline{
{\bf Table I.} Two-dimensional DD algebras.
}
\begin{center}
\begin{tabular}{|c|c|c|}
\hline
     &  $\mathfrak{b}_2\odot R^2$  &  $gl(2)$  
  \\
\hline \hline
 $[x^0,x^1]$ &  $x^1$ &  $x^1$  \\
\hline
 $[X_0,X_1]$ &  $X_0$ &  $\lambda X_1$  \\
\hline
 $[x^0,X_0]$ &  $x^1$ &  $0$  \\
\hline
 $[x^0,X_1]$ &  $-x^0 - X_1$ &  $-X_1$  \\
\hline
 $[x^1,X_0]$ &  $0$ &  $\lambda\,x^1$  \\
\hline
 $[x^1,X_1]$ &  $X_0$ &  $-\lambda x^0 + X_0$  \\
\hline
\end{tabular}
\end{center}
In this case there are only two DD algebras, the ``standard one" (isomorphic to a
semidirect product of the Borel subalgebra $\mathfrak{b}_2$ and $R^2$,
hereafter denoted as
$\mathfrak{b}_2\odot R^2$) and the ``non-standard one" (with $\lambda\neq 0$), which is
isomorphic to $gl(2)$. Let us now describe their quantization (see \cite{HS}
for the connection of both DD algebras with with Poisson-Lie T-duality).

\subsection{The case $\mathfrak{b}_2\odot R^2$ }

A coproduct for this quantum DD algebra with the ``generalized cocommutativity"
property \cite{3dim} would be:
$$ \begin{array}{ll}
 \Delta(x^0)=e^{-z\,x^1}\otimes x^0 + x^0\otimes e^{z\,x^1}, \qquad  
& \Delta(X_0)=1\otimes X_0 + X_0\otimes 1,  \\[0.25cm]
 \Delta(x^1)=1\otimes x^1 + x^1\otimes 1, \qquad  
& \Delta(X_1)= e^{-z\,X_0}\otimes X_1 + X_1\otimes e^{z\,X_0}.   
\nonumber
\end{array}
$$ 
Note that this coassociative coproduct has to be defined in such a way that the dual of
the first-order deformation in the deformation parameter $z$ for the $\{x^0,x^1\}$ (resp.
$\{X_0,X_1\}$) generators provides the dual Lie algebra $\{X_0,X_1\}$ (resp.
$\{x^0,x^1\}$). Explicitly,
$$
\Delta=\Delta_0 + z\,\delta_{DD} + o[z^2],
$$
where $\Delta_0 (y)=1\otimes y + y \otimes 1$ is the primitive coproduct.

The quantum commutation rules compatible with the previous coproduct (and such that they 
deform the Table I relations) can be split into two subsets:

\noindent Quantum subalgebras:
$$
[x^0,x^1]=\frac{\sinh(z x^1)}{z}, \qquad\qquad 
[X_0,X_1]=\frac{\sinh(z X_0)}{z}.
$$ 
Quantum crossed relations:
$$ \begin{array}{ll}
 [x^0,X_0]= \frac{\sinh(z x^1)}{z}, \qquad & [x^1,X_0]= 0,   \\[0.25cm]
 [x^0,X_1]= -\cosh(z\,X_0)\,x^0 - \cosh(z\,x^1)\,X_1, \qquad  
& [x^1,X_1]= \frac{\sinh(z X_0)}{z}.  
\end{array}
$$ 
To our knowledge, this is a new four-generators quantum algebra, whose additional central
element is
$x^1-X_0$. We also remark that, in spite of appearances, also the right hand side of
$[x^0,X_1]$ belongs to the symmetric quantum universal enveloping algebra.

\subsection{The $gl(2)$ case}

The quantization of the remaining DD leads to the following 
coproduct:
$$ \begin{array}{lll}
 \Delta(x^0)=1\otimes x^0 + x^0\otimes 1, \qquad  
& \Delta(X_0)=1\otimes X_0 + X_0\otimes 1,   \\[0.25cm]
 \Delta(x^1)=e^{z\lambda\,\,x^0}\otimes x^1 + x^1\otimes e^{-z\lambda\,x^0}, \qquad  
& \Delta(X_1)=e^{-z\,X_0}\otimes X_1 + X_1\otimes e^{z\,X_0}.   
\nonumber
\end{array}
$$ and quantum commutation rules:

\noindent Subalgebras
$$ [x^0,x^1]=x^1,\qquad \qquad 
[X_0,X_1]=\lambda\,X_1 .
$$ \noindent Crossed relations
$$ \begin{array}{lll}
 [x^0,X_0]= 0, \qquad  & [x^1,X_0]= \lambda x^1,  \\[0.25cm]
 [x^0,X_1]= -X_1, \qquad  & [x^1,X_1]= \frac{\sinh(z (-\lambda\,x^0 + X_0))}{z} .
\end{array}$$ 
We mention that $C=\lambda\,x^0 + X_0$ is a Casimir operator for this quantum algebra,
which has been already studied in \cite{gl2}. Note that the case $\lambda=0$ is
isomorphic to the ``standard case", as it has been pointed out for the DD algebra in
\cite{HS}.

\section{Quantum DD algebras: $\mathfrak{g}=\mathfrak{r}_3(1)$}

In the 3D case we are going to follow the same strategy and notation. For
the sake of simplicity, we shall split Gomez's classification of 3D Lie bialgebras into
different classes according to the structure of the Lie algebra $\mathfrak{g}$.
For each of these classes we shall start by presenting a comprehensive Table of
classical results including the pairs of corresponding 3D dual Lie algebras
$(\mathfrak{g},\mathfrak{g}^\ast)$ (first row), the Lie algebra characterization of
the DD algebras (second row) and Gomez's
notation \cite{gomez} for the all the possible Lie bialgebras (third row). Afterwards, the
non-deformed Lie brackets are presented by columns in the usual basis (following
\cite{gomez}, $\lambda$ is a nonzero essential parameter and $\omega=\pm 1$). 
Next, we shall proceed to the quantization.

We start with the $\mathfrak{g}=\mathfrak{r}_3(1)$ cases, that are explicitly given in
Table II.
We stress that the DDs of type  $(\mathfrak{r}_3(1),\mathfrak{sl}_2)$,
$(\mathfrak{r}_3(1),\mathfrak{so}_3/\mathfrak{sl}_2)$ and 
$(\mathfrak{r}_3(1),\mathfrak{sl}_2)$ are the only DDs for which we have not
succeeded in obtaining a complete quantization, since we have not been
able to construct the quantum coproduct for the $(\mathfrak{g},\delta)$
sub-Lie-bialgebra.  However, we shall see in the sequel that it is possible to identify
these quantum algebras as Drinfel'd twists of previously known
$\mathfrak{sl}_2\oplus \mathfrak{sl}_2$,  
$\mathfrak{so}(1,3)$ and (2+1) Poincar\'e algebra deformations.

\vfill
\eject

\begin{center}
{\bf Table II.} DD algebras with $\mathfrak{g}=\mathfrak{r}_3(1)$. 
\end{center}

\footnotesize
\begin{center}
\begin{tabular}{|c|c|c|c|c|c|c|}
\hline
$(\mathfrak{g},\mathfrak{g}^\ast)$  &  $(\mathfrak{r}_3(1),\mathfrak{sl}_2)$  &  
$(\mathfrak{r}_3(1),\mathfrak{so}_3/\mathfrak{sl}_2)$  & 
$(\mathfrak{r}_3(1),\mathfrak{sl}_2)$   &  $(\mathfrak{r}_3(1),\mathfrak{s}_3(0))$  & 
$(\mathfrak{r}_3(1),\mathfrak{n}_3)$  &  
$(\mathfrak{r}_3(1),\mathfrak{r}_3(-1)/\mathfrak{s}_3(0))$ 
  \\
\hline
DD  &  $\mathfrak{sl}_2\oplus \mathfrak{sl}_2$  &   $\mathfrak{so}(1,3)$ 
& 
$\mathfrak{sl}_2\odot R^3$   &  $\mathfrak{so}(1,3)$  &  $\mathfrak{sl}_2\odot R^3$  &  
$\mathfrak{sl}_2\odot R^3$  
  \\
\hline
 \cite{gomez}  Nr. &  $(1)$  &   $(4)\, \mbox{and} \, (2)$  &   $(3)$
  &  $9$  &  $10$  &  
 $11\, \mbox{and} \, 11'$   
  \\
\hline \hline
 $[x^0,x^1]$ &  $x^1$ &  $x^1$ &  $x^1$ &  $x^1$ &  $x^1$ & $x^1$  \\
\hline
 $[x^0,x^2]$ & $x^2$  &  $x^2$ & $x^2$  & $x^2$  & $x^2$  & $x^2$  \\
\hline
 $[x^1,x^2]$ &  $0$ &  $0$  &   $0$ &  $0$  &  $0$  &  $0$  \\
\hline
 $[X_0,X_1]$ &  $\lambda X_1$  & $\lambda X_2$  & $ X_0$  
&  $\lambda X_2$ & $0$  &  $0$\\
\hline
 $[X_0,X_2]$ &  $-\lambda X_2$  & $-\lambda X_1$  & $X_1$    &   $-\lambda X_1$ &  $X_1$
& $X_1$\\
\hline
 $[X_1,X_2]$ & $ X_0$  &  $\omega X_0$ & $ X_2$   &
 $0$  & $0$  & $ \omega X_0$
\\
\hline
 $[x^0,X_0]$ & $0$  &  $0$ & $ x^1$  &  $0$ &  $0$ &  $0$ \\
\hline
 $[x^0,X_1]$ & $ x^2- X_1$  & $\omega x^2- X_1$  & $- x^0- X_1$   &  $- X_1$  &  $-
X_1$  &  $ \omega x^2- X_1$\\
\hline
 $[x^0,X_2]$ & $- x^1- X_2$  & $-\omega x^1- X_2$   &  $- X_2$   
   & $- X_2$  & $- X_2$  & $- \omega x^1- X_2$\\
\hline
 $[x^1,X_0]$ & $\lambda x^1$  & $-\lambda x^2$  &  $x^2$     & $-\lambda x^2$  & $x^2$  
& $x^2$\\
\hline
 $[x^1,X_1]$ & $-\lambda x^0 + X_0$  & $X_0$  & $X_0$  & $X_0$  & $X_0$  & $X_0$ \\
\hline
 $[x^1,X_2]$ & $0$  & $\lambda x^0$   & $-x^0$       
&  $\lambda x^0$   &  $-x^0$   & $-x^0$ \\
\hline
 $[x^2,X_0]$ &  $-\lambda x^2$  & $\lambda x^1$   &  $0$  & $\lambda x^1$  &  $0$ 
&  $0$\\
\hline
 $[x^2,X_1]$ & $0$  & $-\lambda x^0$  & $ x^2$    & $-\lambda x^0$ & $0$  & $0$ \\
\hline
 $[x^2,X_2]$ & $\lambda x^0 + X_0$  & $X_0$   & $- x^1 +
X_0$   &  $X_0$ &  $X_0$ &  $X_0$ \\
\hline
\end{tabular}
\end{center}

\normalsize

\subsection{Type  $(\mathfrak{r}_3(1),\mathfrak{sl}_2)$ }

In order to characterize the classical $r$-matrix of this DD algebra, we could consider
the following pairing-preserving change of basis:
$$\begin{array}{l}
\hat X_1 = \frac{1}{\sqrt{2}}( X_1 + X_2) + \frac{1}{2\sqrt{2}}( x^1 - x^2),\\[0.25cm]
\hat X_2 =
\frac{1}{\sqrt{2}}( X_1 - X_2) - \frac{1}{2\sqrt{2}}( x^1 + x^2),\\[0.25cm] 
x'^1= \frac{1}{\sqrt{2}}( x^1 + x^2),\\[0.25cm] 
x'^2=\frac{1}{\sqrt{2}}( x^1 -x^2).
\end{array}$$ 
Now, if we recall the
commutation relations of the Lie algebra
$g_{(\m_1,\m_2,\m_3)}$,  a three-parameter family of graded contractions of $so(2,2)$
\cite{so22} (with $\m_i$ being real parameters),  
we
can make the following identification between the DD Lie algebra of type 
$(\mathfrak{r}_3(1),\mathfrak{sl}_2)$  and the 
$g_{(\m_1,\m_2,\m_3)}$ generators:
$$\begin{array}{llll}
& N_3= 2\, X_0, \qquad
& J_3=2\, x^0, \qquad
& J_+=-x'^2 ,\\[0.25cm] 
&N_+=x'^1 ,\qquad
&J_-=-\hat X_1 ,\qquad
&N_-=\hat X_2 .
\end{array}$$ 
It turns out that the DD algebra  $(\mathfrak{r}_3(1),\mathfrak{sl}_2)$  is just a
$so(2,2)$ algebra $g_{(\lambda ,-1/2,\lambda)}$.
Since the associated skew-symmetric classical $r$-matrix is $\tilde r=\frac12
\sum_i{x^i\wedge X_i}$, we get
$$
\tilde r=r_s+ r_{t1}+ r_{t2}=
-\frac12 (N_+\wedge J_- + J_+\wedge N_-) - \frac18 N_3\wedge J_3 +
\,\frac12 N_+\wedge J_+ .
$$
Thus, we have that the quantum DD algebra of type  $(\mathfrak{r}_3(1),\mathfrak{sl}_2)$ 
would be isomorphic to the standard deformation of $so(2,2)$ \cite{so22} generated by
$r_s$ plus two non-commuting twists generated by 
$ r_{t1}+ r_{t2}= - \frac18 N_3\wedge J_3+ \frac12  \,N_+\wedge J_+$.

\subsection{Types $(\mathfrak{r}_3(1),\mathfrak{so}_3/\mathfrak{sl}_2)$}

The two cases are labeled by $\omega\lambda$: when $\omega\lambda>0$ we have the DD
algebra $(\mathfrak{r}_3(1),\mathfrak{so}_3)$, and $\omega\lambda<0$ corresponds to
$(\mathfrak{r}_3(1),\mathfrak{sl}_2)$.

By following the same procedure we
perform the change of basis:
$$
\hat X_1 = X_1 - \frac{\omega}{2}\, x^2,\qquad \hat X_2 = X_2 +\frac{\omega}{2}\, x^1
$$
and we
obtain that the DD algebra $so(3,1)$ is isomorphic to the one with parameters
$g_{(-\lambda,-1/2,\lambda)}$. The associated classical $r$-matrix of this DD algebra is
again
$$
\tilde r=r_s+ r_{t1}+ r_{t2}=
-\frac12 (N_+\wedge J_- + J_+\wedge N_-) - \frac18 N_3\wedge J_3 -
\frac{\omega}{2}\,N_+\wedge J_+  .
$$
Thus, in this case we would have the standard deformation of $\mathfrak{so}(1,3)$ \cite{so22} plus
two non-commuting Drinfel'd twists.

\subsection{Type $(\mathfrak{r}_3(1),\mathfrak{sl}_2)$}

Once again, through the following change of basis
$$
\hat X_1 = X_1+  \, x^0,\qquad \hat X_0 = X_0 -  \, x^1 ,
$$
the Lie algebra parameters turn out to be $g_{(0,-1/2,1)}$, 
which correspond to a $iso(1,2)$ algebra (the 2+1 Poincar\'e algebra in the language
of space-time motion groups). The associated
skew-symmetric classical $r$-matrix is again the standard one \cite{so22} plus a pair of
twists:
$$
\tilde r=r_4+ r_{t6}=
-  \frac12 (N_+\wedge J_- + J_+\wedge N_-) - \frac18 N_3\wedge J_3 -
\frac{1}{2} J_+\wedge J_3 .
$$
We recall that previous works on twist deformations of (3+1) Poincar\'e algebras
were restricted to the study of abelian twists \cite{Luki}.

\subsection{Type $(\mathfrak{r}_3(1),\mathfrak{s}_3(0))$}

Firstly note that this case is the $\omega\to 0$ limit of the type
$(\mathfrak{r}_3(1),\mathfrak{so}_3/\mathfrak{sl}_2)$, and thus it could be considered as
a first step for the quantization of the latter. However, this DD algebra can be directly
quantized and we get:

\noindent Coproduct:
\be\begin{array}{l}
 \Delta(x^0)=1\otimes x^0 + x^0\otimes 1,\\[0.25cm]
 \Delta(x^1)=\cos(z\lambda x^0)\otimes x^1 + 
\sin(z\lambda x^0)\otimes x^2 + x^1\otimes \cos(z\lambda x^0) -
x^2 \otimes \sin(z\lambda x^0),\\[0.25cm]
 \Delta(x^2)=-\sin(z\lambda x^0)\otimes x^1 + 
\cos(z\lambda x^0)\otimes x^2 + x^1\otimes \sin(z\lambda x^0) +
x^2 \otimes \cos(z\lambda x^0),\\[0.25cm]
 \Delta(X_0)=1\otimes X_0 + X_0\otimes 1,\label{tw}\\[0.25cm]
 \Delta(X_1)=e^{zX_0}\otimes X_1 + X_1\otimes e^{-zX_0},\\[0.25cm]
 \Delta(X_2)=e^{zX_0}\otimes X_2 + X_2\otimes e^{-zX_0} .
\end{array}
\ee

\noindent Commutation rules:

\noindent a) Subalgebras:
\be\begin{array}{ll}\label{ctw}
 [x^0,x^1]=x^1, \qquad &[X_0,X_1]=\lambda X_2,\\[0.25cm]
 [x^0,x^2]=x^2,\qquad&[X_0,X_2]=-\lambda X_1,\\[0.25cm]
 [x^1,x^2]=0,\qquad& [X_1,X_2]=0. 
\end{array}
\ee
\noindent b) Crossed relations:
$$ \begin{array}{lll}
 [x^0,X_0]=0,  \quad & [x^1,X_0]=-\lambda x^2, \quad & [x^2,X_0]=\lambda x^1,\\[0.25cm]
 [x^0,X_1]=-X_1,  \quad & [x^1,X_1]=\frac{\sinh(zX_0)}{z} \cos(z\lambda x^0), \quad 
&[x^2,X_1]=-\frac{\sin(z\lambda x^0)}{z} \cosh(z X_0),\\[0.25cm]
 [x^0,X_2]=-X_2,  \quad &  [x^1,X_2]=\frac{\sin(z\lambda x^0)}{z} \cosh(z X_0), \quad 
&[x^2,X_2]=\frac{\sinh(zX_0)}{z} \cos(z\lambda x^0) .
\end{array}
$$ 
The relation with the deformations described in \cite{so22} can be easily derived, and
this DD algebra $\mathfrak{so}(1,3)$ has parameters $g_{(-\lambda,-1/2,\lambda)}$
and canonical $r$-matrix given by
$$
\tilde r=r_s+ r_{t}=-
\frac12 (N_+\wedge J_- + J_+\wedge N_-) - \frac18 N_3\wedge J_3 
$$
Thus, we have again the standard deformation of $\mathfrak{so}(1,3)$ plus a Reshetikhin
twist $r_{t}$.

It is worthy to be emphasized that the canonical basis $\{x^i,X_j\}$ allows a more
manageable description since the
two Hopf subalgebras with classical commutation rules become apparent (compare (\ref{tw}) 
and (\ref{ctw}) with the results in \cite{so22}).

\subsection{Type $(\mathfrak{r}_3(1),\mathfrak{n}_3)$}

This DD algebra is a contraction of both $(\mathfrak{r}_3(1),\mathfrak{sl}_2)$ and
$(\mathfrak{r}_3(1),\mathfrak{r}_3(-1)/\mathfrak{s}_3(0))$. Its quantization is simply
the $\omega
\to 0$ limit of the corresponding expressions in section 4.6.

\noindent Coproduct:
$$ \begin{array}{ll}
 \Delta(x^0)=1\otimes x^0 + x^0\otimes 1, \quad  
& \Delta(X_0)=1\otimes X_0 + X_0\otimes 1,   \\[0.25cm]
 \Delta(x^1)=(1\otimes x^1 + x^1\otimes 1) +
z(x^2\otimes x^0-x^2\otimes x^0), \quad  
& \Delta(X_1)=e^{zX_0}\otimes X_1 + X_1\otimes e^{-zX_0},  \\[0.25cm]
 \Delta(x^2)=1\otimes x^2 + x^2\otimes 1, \quad  
& \Delta(X_2)=e^{zX_0}\otimes X_2 + X_2\otimes
e^{-zX_0}.     
\end{array}
$$ 
\noindent Commutation rules:

\noindent a) Subalgebras:
$$ \begin{array}{ll}
 [x^0,x^1]=x^1, \qquad  & [X_0,X_1]=0,  \\[0.25cm]
 [x^0,x^2]=x^2, \qquad  & [X_0,X_2]=X_1,  \\[0.25cm]
 [x^1,x^2]=0, \qquad  & [X_1,X_2]=0 .   
\nonumber
\end{array}
$$ 
\noindent b) Crossed relations:
$$ \begin{array}{lll}
 [x^0,X_0]=0, \qquad  & [x^1,X_0]= x^2, \qquad & [x^2,X_0]=0,\\[0.25cm]
 [x^0,X_1]=-X_1, \qquad  & [x^1,X_1]=\frac{\sinh(zX_0)}{z}, \qquad 
&[x^2,X_1]=0,\\[0.25cm]
 [x^0,X_2]=-X_2, \qquad  & [x^1,X_2]=-x^0 \cosh{z X_0},  \qquad
 & [x^2,X_2]=\frac{\sinh(zX_0)}{z} .
\end{array}
$$ 
Once again,
we get the characterization $g_{(0,-1/2,1)}$.
Not that by using the transformation ${X_1}'=-\lambda X_1$, ${x^1}'=- \frac 1 {\lambda}
x^1$, the type $(\mathfrak{r}_3(1),\mathfrak{n}_3)$ is just the limit for $\lambda \to 0$
of the known type
$(\mathfrak{r}_3(1),\mathfrak{s}_3(0))$. The associated
skew-symmetric classical $r$-matrix is:
$$
\tilde r=r_s+ r_{t}=
- \frac12 (N_+\wedge J_- + J_+\wedge N_-) -  \frac18 N_3\wedge J_3 .
$$
Thus, we have again the Poincar\'e
analogue to the previous case: the known standard deformation \cite{so22} plus a twist.

\subsection{Types $(\mathfrak{r}_3(1),\mathfrak{r}_3(-1)/\mathfrak{s}_3(0))$}

We remark that we have $(\mathfrak{r}_3(1),\mathfrak{r}_3(-1))$ for $\omega=+1$, and 
$(\mathfrak{r}_3(1),\mathfrak{s}_3(0))$ for $\omega=-1$.
The full quantization is as follows.

\noindent Coproduct:
\[\begin{array}{l} 
 \Delta(x^0)=
   \cosh( \sqrt{ \omega} z x^2)\otimes x^0 +
    \sqrt{ \omega} \sinh( \sqrt{ \omega} z x^2)\otimes x^1 \\[0.15cm]
\hskip2cm  + x^0 \otimes\cosh( \sqrt{ \omega} z x^2) -
    \sqrt{ \omega} x^1\otimes \sinh( \sqrt{ \omega} z x^2), \\[0.25cm]
 \Delta(x^1)=
  \frac{\sinh( \sqrt{ \omega} z x^2)}{ \sqrt{ \omega}} \otimes x^0 +    
  \cosh( \sqrt{ \omega} z x^2)\otimes x^1  \\[0.15cm]
\hskip2cm   - x^0\otimes \frac{\sinh( \sqrt{ \omega} z x^2)}{ \sqrt{ \omega}} +
  x^1 \otimes\cosh( \sqrt{ \omega} z x^2),\\[0.25cm]
\Delta(x^2)=1\otimes x^2 + x^2\otimes 1,\\[0.25cm]
 \Delta(X_0)=1\otimes X_0 + X_0\otimes 1,\\[0.25cm]
\Delta(X_1)=e^{zX_0}\otimes X_1 + X_1\otimes e^{-zX_0},\\[0.25cm]
 \Delta(X_2)=e^{zX_0}\otimes X_2 + X_2\otimes e^{-zX_0}.
\end{array}\]

\noindent Commutation rules:

\noindent a) Subalgebras:
$$ \begin{array}{lll}
 [x^0,x^1]=x^1 \,\cosh{[ \sqrt{ \omega} z x^2]}, \qquad  & [X_0,X_1]=0, \\[0.25cm] 
 [x^0,x^2]= \frac {\sinh {[ \sqrt{ \omega} z x^2]}} { \sqrt{ \omega} z}, \qquad  
& [X_0,X_2]=X_1,  \\[0.25cm]
 [x^1,x^2]= 0 , \qquad  & [X_1,X_2]=\omega \frac {\sinh(2zX_0)} {2z}.    
\end{array}
$$ 
\noindent b) Crossed relations:
$$\begin{array}{ll}
 [x^0,X_0]=0,& \\[0.25cm] 
 [x^0,X_1]=\frac { \sqrt{ \omega}} {z} \sinh{ \left[ \sqrt{ \omega}\, z
x^2\right]}\,\cosh{[z X_0]} - \cosh{\left[  \sqrt{ \omega}\, z x^2\right]}
X_1,& \\[0.25cm]  
 [x^0,X_2]=-\omega x^1 \cosh{ (z X_0)} -
X_2 \cosh{( \sqrt{ \omega} z x^2)},& \\[0.25cm] 
 [x^1,X_0]= \frac{\sinh{\left[ \sqrt{ \omega} z x^2\right]}}
{ \sqrt{ \omega} \,z}, \qquad  & [x^2,X_0]=0,   \\[0.25cm] 
 [x^1,X_1]= \cosh{\left[  \sqrt{ \omega}\, z x^2\right]}\,\,
\frac {\sinh{[z X_0]}} {z}, \qquad  & [x^2,X_1]=0,     \\[0.25cm] 
 [x^1,X_2]=-x^0 \cosh{(z X_0)}, \qquad  & [x^2,X_2]=\frac{\sinh(z X_0)}{z}.     
\end{array}
$$

Note that the use of a fully symmetric basis is implicit since
\[\begin{array}{l}
 [x^0,X_2]=- \omega x^1 \cosh{ (z X_0)} -
X_2 \cosh{( \sqrt{ \omega} z x^2)} \\[0.25cm] 
\hskip1.45cm =  - \omega \,\mbox{Sym}[x^1 \cosh{ (z X_0)}] -
\mbox{Sym}[X_2 \cosh{( \sqrt{ \omega} z x^2)}],
\end{array}\]
where ${\rm Sym}$ is a linear operator such that
$$
{\rm Sym}\;\{A_1\dots A_n\} := \frac 1{n!}\sum_{p\in {\rm S_n}} 
p (A_1\dots A_n),
$$
with ${\rm S_n}$ the permutation group of $n$ elements (see \cite{3dim}).

In this case the DD algebra is a $iso(2,1)$ algebra with parameters $g_{(1,-1/2,0)}$
and the classical $r$-matrix would be
$$
\tilde r=r_s+ r_{t1}+ r_{t2}=
-\frac12 (N_+\wedge J_- + J_+\wedge N_-) - \frac18 N_3\wedge J_3 -
 \frac{\omega}{2}\,N_+\wedge J_+ .
$$
Once again, we have the standard deformation of the  2+1 Poincar\'e algebra plus two
non-commuting twists, with the same interpretation as for the type 
$(\mathfrak{r}_3(1),\mathfrak{sl}_2)$ for
$\mathfrak{so}(1,3)$.

\section{Quantum DD algebras: $\mathfrak{g}=\mathfrak{r}_3(\rho)$}

The second set of quantum DD algebras corresponds to the quantization of those DD
algebras with $\mathfrak{g}=\mathfrak{r}_3(\rho)$. The three possible cases are described
in Table III and can be fully quantized.

\bigskip
 
\centerline{
{\bf Table II.} DD algebras with $\mathfrak{g}=\mathfrak{r}_3(\rho)$.
}

\footnotesize
\begin{center}
\begin{tabular}{|c|c|c|c|}
\hline
$(\mathfrak{g},\mathfrak{g}^\ast)$  &  $(\mathfrak{r}_3(\rho),\mathfrak{n}_3)$  &  
$(\mathfrak{r}_3(\rho),\mathfrak{r}_3(-\rho))$  & 
$(\mathfrak{r}_3(\rho),\mathfrak{r}_3(-\rho))$  
  \\
\hline
DD &  $\mathfrak{r}_3(\rho)\odot R^3$   &   $\mathfrak{r}_3(\rho)\odot R^3$
& 
$sl_2\oplus sl_2$    \\
\hline
 \cite{gomez}  Nr. &  $5$  &   $6$  &   $7$   \\
\hline \hline
 $[x^0,x^1]$ &  $x^1$ &  $x^1$ &  $x^1$\\
\hline
 $[x^0,x^2]$ & $\rho\,x^2$  &  $\rho\,x^2$ & $\rho\,x^2$ \\
\hline
 $[x^1,x^2]$ &  $0$ &  $0$  &   $0$ \\
\hline
 $[X_0,X_1]$ &  $0$  & $ X_0$  & $\lambda X_1$  \\
\hline
 $[X_0,X_2]$ &  $0$  & $0$  &  $-\lambda \rho X_2$ \\
\hline
 $[X_1,X_2]$ & $ (1+\rho) X_0$  &  $\rho X_2$ & $0$  \\
\hline
 $[x^0,X_0]$ & $0$  &  $x^1$ & $0$ \\
\hline
 $[x^0,X_1]$ & $(1+\rho) x^2 - X_1$  & $- x^0 - X_1$  & $- X_1$  \\
\hline
 $[x^0,X_2]$ & $-(1+\rho) x^1 - \rho X_2$  & $- \rho X_2$   &  $- \rho X_2$ \\
\hline
 $[x^1,X_0]$ & $0$  & $0$  &  $\lambda x^1$  \\
\hline
 $[x^1,X_1]$ & $X_0$  & $X_0$  & $-\lambda x^0 + X_0$  \\
\hline
 $[x^1,X_2]$ & $0$  & $0$   & $0$     \\
\hline
 $[x^2,X_0]$ &  $0$  & $0$   &  $-\lambda \rho x^2$  \\
\hline
 $[x^2,X_1]$ & $0$  & $\rho \,x^2$  & $0$   \\
\hline
 $[x^2,X_2]$ & $\rho\,X_0$  & $\rho (-x^1 + X_0)$   & $\rho(\lambda x^0 + X_0)$  \\
\hline
\end{tabular}
\end{center}

\normalsize
\subsection{Case $(\mathfrak{r}_3(\rho),\mathfrak{n}_3)$}

Coproduct:
$$ \begin{array}{l}
 \Delta(x^0)=1\otimes x^0 + x^0\otimes 1 - z(1+\rho)(x^1\otimes x^2-x^2\otimes x^1), 
\\[0.25cm]
\Delta(X_0)=1\otimes X_0 + X_0\otimes 1,   \\[0.25cm]
 \Delta(x^1)=1\otimes x^1 + x^1\otimes 1, \\[0.25cm] 
\Delta(X_1)=e^{zX_0}\otimes X_1 + X_1\otimes e^{-zX_0},  \\[0.25cm]
 \Delta(x^2)=1\otimes x^2 + x^2\otimes 1, \\[0.25cm] 
 \Delta(X_2)=e^{z\rho X_0}\otimes X_2 + X_2\otimes
e^{-z\rho X_0}.     \nonumber
\end{array}
$$ 

\noindent Commutation rules:

\noindent a) Subalgebras:
$$ \begin{array}{lll}
 [x^0,x^1]=x^1, \qquad  & [X_0,X_1]=0, \\[0.25cm]
 [x^0,x^2]=\rho x^2, \qquad  & [X_0,X_2]=0, \\[0.25cm]
 [x^1,x^2]=0, \qquad  & [X_1,X_2]=\frac{\sinh(z(1+\rho)X_0)}{z}.    
\end{array}
$$ 
\noindent b) Crossed relations:
$$ \begin{array}{lll}
 [x^0,X_0]=0, \quad  & [x^1,X_0]= 0, \quad & [x^2,X_0]=0,\\[0.25cm]
 [x^0,X_1]=(1+\rho)x^2 \cosh(z X_0) - X_1, \quad  
& [x^1,X_1]=\frac{\sinh(z X_0)}{z}, \quad & [x^2,X_1]=0,\\[0.25cm]
 [x^0,X_2]=-(1+\rho)x^1 \cosh(z \rho X_0) - \rho X_2, \quad  & [x^1,X_2]=0,  \quad 
& [x^2,X_2]=\frac{\sinh(z\rho X_0)}{z}.
\end{array}
$$ 

\subsection{Case $(\mathfrak{r}_3(\rho),\mathfrak{r}_3(-\rho))$}

Coproduct:
$$ \begin{array}{ll}
 \Delta(x^0)=e^{zx^1}\otimes x^0 + x^0\otimes e^{-zx^1}, \qquad  
& \Delta(X_0)=1\otimes X_0 + X_0\otimes 1,  \\[0.25cm]
 \Delta(x^1)=1\otimes x^1 + x^1\otimes 1, \qquad  
& \Delta(X_1)=e^{zX_0}\otimes X_1 + X_1\otimes e^{-zX_0},   \\[0.25cm]
 \Delta(x^2)=e^{-z\rho x^1}\otimes x^2 + x^2\otimes e^{z\rho x^1}, \qquad 
 & \Delta(X_2)=e^{z\rho X_0}\otimes X_2 + X_2\otimes e^{-z\rho X_0}.     
\end{array}
$$ 
\noindent Commutation rules:

\noindent a) Subalgebras:
$$ \begin{array}{ll}
 [x^0,x^1]=\frac{\sinh(z x^1)}{z}, \qquad  & [X_0,X_1]=\frac{\sinh(z X_0)}{z},  \\[0.25cm]
 [x^0,x^2]=\rho x^2 \cosh(z x^1), \qquad  & [X_0,X_2]=0,  \\[0.25cm]
 [x^1,x^2]=0, \qquad  & [X_1,X_2]=\rho X_2 \cosh(z X_0).     \nonumber
\end{array}
$$ 
\noindent b) Crossed relations:
$$
\begin{array}{lll}
 [x^0,X_0]=\frac{\sinh(z x^1)}{z}, \quad  & [x^1,X_0]= 0, \quad & [x^2,X_0]=0,\\[0.25cm]
 [x^0,X_1]= -\cosh(z\,X_0)\,x^0 \quad  
& [x^1,X_1]=\frac{\sinh(z X_0)}{z}, \quad & [x^2,X_1]=\rho x^2 \cosh(z X_0),\\[0.1cm]
\hskip2cm  - \cosh(z\,x^1)\,X_1,&& \\[0.25cm]
 [x^0,X_2]=- \rho X_2 \cosh(z x^1), \quad  &
[x^1,X_2]=0, 
\quad  &  [x^2,X_2]=-\frac{\sinh(z \rho (-x^1+X_0))}{z}. 
\end{array}
$$ 
The symmetrization prescription is again preserved despite the non-symmetric shape of
some brackets. Note that this DD algebra is self-dual for
$\rho=0$.

\subsection{Case $(\mathfrak{r}_3(\rho),\mathfrak{r}_3(-\rho))$}

Coproduct:
$$ \begin{array}{ll}
 \Delta(x^0)=1\otimes x^0 + x^0\otimes 1, \qquad  
& \Delta(X_0)=1\otimes X_0 + X_0\otimes 1, \\[0.25cm]
 \Delta(x^1)=e^{-z\lambda x^0}\otimes x^1 + x^1\otimes e^{z\lambda x^0}, \qquad  
&\Delta(X_1)=e^{zX_0}\otimes X_1 + X_1\otimes e^{-zX_0},   \\[0.25cm]
 \Delta(x^2)=e^{z\lambda\rho  x^0}\otimes x^2 + x^2\otimes e^{-z\lambda\rho x^0}, \qquad  
& \Delta(X_2)=e^{z\rho X_0}\otimes X_2 + X_2\otimes e^{-z\rho X_0}.     
\end{array}
$$ 
\noindent Commutation rules:

\noindent a) Subalgebras:
$$ \begin{array}{ll}
 [x^0,x^1]=x^1, \qquad  & [X_0,X_1]=\lambda X_1,  \\[0.25cm]
 [x^0,x^2]=\rho x^2, \qquad  & [X_0,X_2]=-\lambda \rho X_2,   \\[0.25cm]
 [x^1,x^2]=0, \qquad  & [X_1,X_2]=0     
\end{array}
$$ 
\noindent b) Crossed relations:
$$ \begin{array}{lll}
 [x^0,X_0]=0, \qquad  & [x^1,X_0]= \lambda x^1, \qquad 
& [x^2,X_0]=-\lambda\rho x^2,\\[0.25cm]
 [x^0,X_1]=- X_1, \qquad  & [x^1,X_1]=\frac{\sinh(z (-\lambda x^0 + X_0))}{z}, \qquad 
& [x^2,X_1]=0,\\[0.25cm]
 [x^0,X_2]=- \rho X_2, \qquad  & [x^1,X_2]=0,  \qquad 
&  [x^2,X_2]=\frac{\sinh(z\rho (\lambda x^0 + X_0))}{z}.
\end{array}
$$ 
Note that the quantum DD algebra corresponding to $gl(2)$  is included as a subalgebra
in several different ways. Again, this DD algebra is self-dual for
$\rho=0$.

\section{Quantum DD algebras:
$\mathfrak{g}=\{\mathfrak{r}_3(-1),\mathfrak{r}_3'(1)\}$}

The following set of quantum DD algebras (see Table IV below) can be also completely
quantized by using a direct quantization approach. Note that all the DD algebras
have a semidirect product structure.

\vfill\eject

\centerline{
{\bf Table IV.} DD algebras with $\mathfrak{g}=\{\mathfrak{r}_3(-1),\mathfrak{r}_3'(1)\}$.
}
\footnotesize
\begin{center}
\begin{tabular}{|c|c|c|c|c|c|}
\hline
$(\mathfrak{g},\mathfrak{g}^\ast)$  &  $(\mathfrak{r}_3(-1),\mathfrak{n}_3)$  &  
$(\mathfrak{r}_3(-1),\mathfrak{r}_3'(1))$  & 
$(\mathfrak{r}_3(-1),\mathfrak{r}_3'(1))$   &  
$(\mathfrak{r}_3'(1),\mathfrak{n}_3)$  & 
$(\mathfrak{r}_3'(1),\mathfrak{n}_3)$   
  \\
\hline
DD  &  $\mathfrak{r}_3'(1)\odot R^3$  &   $\mathfrak{r}_3'(1)\odot R^3$ 
& 
$\mathfrak{sl}_2\odot R^3$ &  $\mathfrak{r}_3'(1)\odot R^3$  &   $sl_2\odot R^3$  \\
\hline
 \cite{gomez} Nr. &  $5'$  &   $8$  &   $(14)$   &  $12$  &   $13$ \\
\hline \hline
 $[x^0,x^1]$ &  $x^1$ &  $x^1$ &  $x^1$ &  $x^1$ &  $x^1$\\
\hline
 $[x^0,x^2]$ & $- x^2$  &  $- x^2$ & $- x^2$& $x^1 + x^2$  &  $x^1 + x^2$\\
\hline
 $[x^1,x^2]$ &  $0$ &  $0$  &   $0$ &  $0$ &  $0$  \\
\hline
 $[X_0,X_1]$ &  $0$  & $   X_0$  & $  X_0$ &  $0$ &  $\lambda X_2$  \\
\hline
 $[X_0,X_2]$ &  $0$  & $0$  &  $-  \lambda X_0$ &  $0$  & $0$ \\
\hline
 $[X_1,X_2]$ & $ X_0$  &  $X_0 -   X_2$ & $X_0 -  \lambda X_1 -   X_2$   & $\omega X_0$ 
&  $0$
\\
\hline
 $[x^0,X_0]$ & $0$  &  $  x^1$ & $  x^1 -  \lambda x^2$ & $0$  &  $0$\\
\hline
 $[x^0,X_1]$ & $x^2 - X_1$  & $- X_1 -   x^0 + x^2$  & $- X_1 -   x^0 + x^2$ & $\omega
x^2 - X_1 - X_2$  & $- X_1 - X_2$  \\
\hline
 $[x^0,X_2]$ & $- x^1 +  X_2$  & $- x^1 + X_2$   &  $ \lambda x^0 - x^1 + X_2$ & $-
\omega x^1 -  X_2$  &
$- X_2$  \\
\hline
 $[x^1,X_0]$ & $0$  & $0$  &  $0$ & $0$  & $0$ \\
\hline
 $[x^1,X_1]$ & $X_0$  & $X_0$  & $- \lambda x^2 + X_0$  & $X_0$  & $X_0$ \\
\hline
 $[x^1,X_2]$ & $0$  & $0$   & $ \lambda x^1$  & $0$  & $0$    \\
\hline
 $[x^2,X_0]$ &  $0$  & $0$   &  $0$ &  $0$  & $\lambda x^1$ \\
\hline
 $[x^2,X_1]$ & $0$  & $-  x^2$  & $-  x^2$ & $X_0$  & $X_0 - \lambda x^0$ \\
\hline
 $[x^2,X_2]$ & $-X_0$  & $  x^1 - X_0$  & $  x^1 - X_0$ & $X_0$  & $X_0$ \\
\hline
\end{tabular}
\end{center}

\normalsize

\subsection{Case $(\mathfrak{r}_3(-1),\mathfrak{n}_3)$}

Coproduct:
$$ \begin{array}{ll}
 \Delta(x^0)=1\otimes x^0 + x^0\otimes 1 - z (x^1\otimes x^2 - x^2\otimes x^1), \qquad  
&\Delta(X_0)=1\otimes X_0 + X_0\otimes 1, \\[0.25cm]
 \Delta(x^1)=1\otimes x^1 + x^1\otimes 1, \qquad  
& \Delta(X_1)=e^{zX_0}\otimes X_1 + X_1\otimes e^{-zX_0},  \\[0.25cm]
 \Delta(x^2)=1\otimes x^2 + x^2\otimes 1, \qquad  
& \Delta(X_2)=e^{- z X_0}\otimes X_2 + X_2\otimes e^{z X_0}.     
\end{array}
$$ 
\noindent Commutation rules:

\noindent a) Subalgebras:
$$ \begin{array}{ll}
 [x^0,x^1]=x^1, \qquad  & [X_0,X_1]=0,  \\[0.25cm]
 [x^0,x^2]=- x^2, \qquad  & [X_0,X_2]=0,  \\[0.25cm]
 [x^1,x^2]=0, \qquad  & [X_1,X_2]=X_0.     
\end{array}
$$ 
\noindent b) Crossed relations:
$$ \begin{array}{lll}
 [x^0,X_0]=0, \qquad  & [x^1,X_0]=0, \qquad & [x^2,X_0]=0,\\[0.25cm]
 [x^0,X_1]=x^2 \cosh(z X_0) - X_1, \qquad  
& [x^1,X_1]=\frac{\sinh(z X_0 )}{z}, \qquad & [x^2,X_1]=0,\\[0.25cm]
 [x^0,X_2]=-x^1 \cosh(z X_0) + X_2, \qquad  
& [x^1,X_2]=0,  \qquad &  [x^2,X_2]=-\frac{\sinh(z X_0)}{z}.
\end{array}
$$ 
\subsection{Case $(\mathfrak{r}_3(-1),\mathfrak{r}_3'(1))$}

Coproduct:
\[\begin{array}{l}
 \Delta(x^0)=e^{  z x^1}\otimes x^0 + x^0\otimes e^{-  z x^1} - z (x^1
e^{  z x^1}\otimes x^2 - x^2\otimes x^1 e^{-  z x^1}) ,\\[0.25cm]
\Delta(x^1)=1\otimes
x^1 + x^1\otimes 1,\\[0.25cm]
 \Delta(x^2)=e^{  z x^1}\otimes x^2 + x^2\otimes e^{-  z x^1},\\[0.25cm]
 \Delta(X_0)=1\otimes X_0 + X_0\otimes 1,\\[0.25cm]
 \Delta(X_1)=e^{zX_0}\otimes X_1 + X_1\otimes e^{-zX_0},\\[0.25cm]
 \Delta(X_2)=e^{- z X_0}\otimes X_2 + X_2\otimes e^{z X_0} .
\end{array}\]

\noindent Commutation rules:

\noindent a) Subalgebras:
$$ \begin{array}{ll}
 [x^0,x^1]=\frac{\sinh(  z x^1 )}{  z}, \qquad  & [X_0,X_1]= 
\frac{\sinh(z X_0)}{z},  \\[0.25cm]
 [x^0,x^2]=- x^2 \cosh(  z x^1), \qquad  & [X_0,X_2]=0,  \\[0.25cm]
 [x^1,x^2]=0, \qquad  & [X_1,X_2]=X_0 -   X_2 \cosh(z X_0).     
\end{array}
$$ 
\noindent b) Crossed relations:
\[\begin{array}{ll} [x^0,X_0]=\frac{\sinh(  z x^1 )}{z},&\\[0.25cm]
 [x^0,X_1]=- X_1 \cosh(  z x^1) -   x^0 \cosh(z X_0) + x^2 \cosh(z X_0),&\\[0.25cm]
 [x^0,X_2]=- x^1 \cosh(z(  x^1 - X_0)) + X_2 \cosh(  z x^1),&\\[0.25cm]
 [x^1,X_0]=0, \qquad  & [x^2,X_0]=0,  \\[0.25cm]
 [x^1,X_1]=\frac{\sinh(z X_0 )}{z}, \qquad  & [x^2,X_1]=-  x^2 \cosh(z X_0),  \\[0.25cm]
 [x^1,X_2]=0, \qquad  & [x^2,X_2]=\frac{\sinh(z (  x^1 - X_0))}{z}.     
\end{array}\]
\subsection{Case $(\mathfrak{r}_3(-1),\mathfrak{r}_3'(1))$}

In order to get the explicit quantization, the following pairing-preserving change of
basis turns out to be convenient:
$$ \begin{array}{ll}
 Y_0=X_0, \qquad  & y^0= x^0,  \\[0.25cm]
 Y_1 = \frac{1}{\sqrt{1+\lambda^2}}( \lambda X_1 + X_2) , \qquad  
& y^1=  \frac{1}{\sqrt{1+\lambda^2}}( \lambda\,x^1 + x^2),   \\[0.25cm]
 Y_2 = \frac{1}{\sqrt{1+\lambda^2}}( X_1 - \lambda  X_2) , \qquad  
& y^2= \frac{1}{\sqrt{1+\lambda^2}}( x^1 - \lambda \, x^2).     
\end{array}
$$ 
In this new basis, the quantum DD algebra reads:

\noindent Coproduct:
\[\begin{array}{l}
 \Delta(y^0)=e^{-  z \sqrt{1+\lambda^2} y^2}\otimes y^0 + y^0\otimes e^{  z
\sqrt{1+\lambda^2} y^2} - z(y^1\otimes y^2\,e^{  z
\sqrt{1+\lambda^2} y^2} - y^2\,e^{-  z \sqrt{1+\lambda^2} y^2}\otimes y^1) \\[0.25cm]
\Delta(y^1)=e^{-  z \sqrt{1+\lambda^2} y^2}\otimes y^1 + y^1\otimes e^{  z
\sqrt{1+\lambda^2} y^2},\\[0.25cm]
\Delta(y^2)=1\otimes y^2 + y^2\otimes 1,\\[0.25cm]
 \Delta(Y_0)=1\otimes Y_0 + Y_0\otimes 1,v
 \Delta(Y_1)= \frac{1}{{1+\lambda^2}} \big\{ ( e^{ z Y_0} 
+ \lambda^2 e^{ -z Y_0})\otimes Y_1
+ Y_1\otimes ( e^{ -z Y_0} + \lambda^2 e^{  z Y_0}) v \\[0.15cm]
\hskip2cm - 2\lambda \sinh(z\,Y_0)\otimes Y_2 + 2\lambda Y_2
\otimes\sinh(z\,Y_0)\big\},\\[0.25cm]
\Delta(Y_2)= \frac{1}{{1+\lambda^2}} \big\{ ( e^{ -z Y_0} + \lambda^2 e^{ z Y_0})\otimes
Y_2 + Y_2\otimes ( e^{ z Y_0} + \lambda^2 e^{ - z Y_0}) \\[0.15cm] 
\hskip2cm 
- 2\lambda \sinh(z\,Y_0)\otimes Y_1 + 2\lambda Y_1 \otimes\sinh(z\,Y_0)\big\}.
\end{array}\]

\noindent Commutation rules:

\noindent a) Subalgebras:
\[\begin{array}{l}
 [y^0,y^1]=- \frac{1-\lambda^2}{1+\lambda^2} \, y^1 \cosh( \,z\,\sqrt{1+\lambda^2}\,y^2)
+ \frac{2\lambda}{1+\lambda^2} \frac{\sinh(2  \,z\,\sqrt{1+\lambda^2}\,y^2)}{2
 \,z\,\sqrt{1+\lambda^2}},\\[0.25cm]
 [y^0,y^2]=\frac{2\lambda}{1+\lambda^2} y^1 +  \frac{1-\lambda^2}{1+\lambda^2}
\frac{\sinh( \,z\,\sqrt{1+\lambda^2}\,y^2)}{
 \,z\,\sqrt{1+\lambda^2}},\\[0.25cm]
 [y^1,y^2]=0,\\[0.25cm]
 [Y_0,Y_1]=0,\\[0.25cm]
 [Y_0,Y_2]=  \,\sqrt{1+\lambda^2}\,\frac{\sinh(z Y_0)}{z},\\[0.25cm]
 [Y_1,Y_2]=-Y_0 +   \,\sqrt{1+\lambda^2} \,Y_1 \,\cosh(z Y_0).\nonumber
\end{array}\]

\noindent b) Crossed relations:
\[\begin{array}{l}
 [y^0,Y_0]=\frac{\sinh( \,z\,\sqrt{1+\lambda^2}\,y^2)}{z},\\[0.25cm]
 [y^0,Y_1]=\left( \frac{1-\lambda^2}{1+\lambda^2} \,Y_1 -
\frac{2\lambda}{1+\lambda^2}\,Y_2\right) \cosh( \,z\,\sqrt{1+\lambda^2}\,y^2)  \\[0.15cm]
 \hskip2cm -  y^2\,\cosh( \,z\,\sqrt{1+\lambda^2}\,y^2) +
\frac{1-\lambda^2}{1+\lambda^2}\,y^2\,\sinh( \,z\,\sqrt{1+\lambda^2}\,y^2)\,\sinh(z
Y_0),\\[0.25cm]
 [y^0,Y_2]=- \,\sqrt{1+\lambda^2}\,y^0\,\cosh(z Y_0) -
\left(\frac{2\lambda}{1+\lambda^2}\,Y_1 - 
\frac{1-\lambda^2}{1+\lambda^2} \,Y_2\right)\cosh( \,z\,\sqrt{1+\lambda^2}\,y^2)  
\\[0.15cm]
\hskip2cm - \frac{2\lambda}{1+\lambda^2}\, y^2\,\sinh( \,z\,\sqrt{1+\lambda^2}\,y^2)+ 
y^1\,\cosh(z Y_0),\\[0.25cm]
 [y^1,Y_0]=0,\\[0.25cm]
 [y^1,Y_1]=\frac{\sinh( \,z\,\sqrt{1+\lambda^2}\,y^2)}{
z}\,\cosh(z Y_0) - \frac{1-\lambda^2}{1+\lambda^2}\,\frac{\sinh(z Y_0)}{
z}\,\cosh( \,z\,\sqrt{1+\lambda^2}\,y^2)    ,\\[0.25cm]
 [y^1,Y_2]=-  \,\sqrt{1+\lambda^2}\,y^1\,\cosh(z Y_0) +
\frac{2\lambda}{1+\lambda^2}\,\frac{\sinh(z Y_0)}{ z}\,
\cosh( \,z\,\sqrt{1+\lambda^2}\,y^2)   ,\\[0.25cm]
 [y^2,Y_1]=\frac{2\lambda}{1+\lambda^2}\,\frac{\sinh(z Y_0)}{ z},\\[0.25cm] 
 [y^2,Y_2]=\frac{1-\lambda^2}{1+\lambda^2}\,\frac{\sinh(z Y_0)}{ z}.
\end{array}\]

\subsection{Case $(\mathfrak{r}_3'(1),\mathfrak{n}_3)$}

Coproduct:
\[\begin{array}{l}
 \Delta(x^0)=1\otimes x^0 + x^0\otimes 1 + \omega z (x^1\otimes x^2
 - x^2\otimes x^1),\\[0.25cm]
 \Delta(x^1)=1\otimes x^1 + x^1\otimes 1,\\[0.25cm]
 \Delta(x^2)=1\otimes x^2 + x^2\otimes 1,\\[0.25cm]
 \Delta(X_0)=1\otimes X_0 + X_0\otimes 1,\\[0.25cm]
 \Delta(X_1)=e^{-zX_0}\otimes X_1 + X_1\otimes e^{zX_0} - z(X_0\,e^{-zX_0}\otimes X_2 - 
X_2\otimes  X_0\,e^{zX_0}) ,\\[0.25cm]
 \Delta(X_2)=e^{-z X_0}\otimes X_2 + X_2\otimes e^{z X_0},
\end{array}\]

\noindent Commutation rules:

\noindent a) Subalgebras:
$$ \begin{array}{ll}
 [x^0,x^1]=x^1, \qquad  & [X_0,X_1]=0,  \\[0.25cm]
 [x^0,x^2]=x^1+  x^2, \qquad  & [X_0,X_2]=0,  \\[0.25cm]
 [x^1,x^2]=0, \qquad  & [X_1,X_2]=\omega \frac{\sinh(2 z X_0)}{2 z}.     
\end{array}
$$ 

\noindent b) Crossed relations:
\[\begin{array}{ll}
 [x^0,X_0]=0,&\\[0.25cm]
 [x^0,X_1]=- X_1 -X_2 + \omega x^2\cosh(z X_0)&\\[0.15cm] 
\hskip2.3cm - \omega z\,X_0\, x^1\sinh(z X_0),\\[0.25cm] 
 [x^0,X_2]=- X_2 - \omega x^1\cosh(z X_0) ,&\\[0.25cm] 
[x^1,X_0]= 0, \qquad   & [x^2,X_0]=0,  \\[0.25cm]
 [x^1,X_1]=\frac{\sinh(z X_0)}{z}, \qquad  & [x^2,X_2]=\frac{\sinh(z X_0)}{z}.     
\end{array}\]

\subsection{Case $(\mathfrak{r}_3'(1),\mathfrak{n}_3)$}

Coproduct:
\[\begin{array}{l}
 \Delta(x^0)=1\otimes x^0 + x^0\otimes 1  ,\\[0.25cm] 
\Delta(x^1)=1\otimes x^1 + x^1\otimes 1 ,\\[0.25cm]  
 \Delta(x^2)=1\otimes x^2 + x^2\otimes 1 + \lambda z\; (x^0\otimes x^1-  
x^1 \otimes x^0),\\[0.25cm]  
\Delta(X_0)=1\otimes X_0 + X_0\otimes 1,\\[0.25cm] 
 \Delta(X_1)=e^{-\,z\,X_0}\otimes X_1 + X_1\otimes 
 e^{z\,X_0} - z\; (X_0 e^{-\,z\,X_0}\otimes X_2-  X_2 \otimes X_0 e^{z\,X_0}),\\[0.25cm] 
\Delta(X_2)= e^{-\,z\,X_0}\otimes X_2 + X_2\otimes   e^{z\,X_0} .
\end{array}\]

\noindent Commutation rules:

\noindent a) Subalgebras:
$$ \begin{array}{ll}
 [x^0,x^1]=x^1, \qquad  & [X_0,X_1]=\lambda X_2,  \\[0.25cm] 
 [x^0,x^2]=x^1+  x^2, \qquad  & [X_0,X_2]=0,   \\[0.25cm] 
 [x^1,x^2]=0, \qquad  & [X_1,X_2]=0.     
\end{array}
$$ 

\noindent b) Crossed relations:
$$ \begin{array}{lll}
 [x^0,X_0]=0, \qquad  & [x^1,X_0]= 0, \qquad & [x^2,X_0]=\lambda x^1,\\[0.25cm] 
 [x^0,X_1]=- X_1 -X_2  , \qquad  & [x^1,X_1]=\frac{\sinh(z X_0)}{z}, \qquad 
& [x^2,X_1]= (X_0 - \lambda x^0)\,\cosh(z X_0),\\[0.25cm] 
 [x^0,X_2]=- X_2 , \qquad  & [x^1,X_2]=0,  \qquad 
&  [x^2,X_2]=\frac{\sinh(z X_0)}{z}.  
\end{array}
$$ 

\section{Quantum DD algebras:
$\mathfrak{g}=\{\mathfrak{s}_3(\mu),\mathfrak{s}_3(0),\mathfrak{n}_3\}$}

Finally, we shall consider the four remaining cases of DD algebras summarized in Table V.
All of them can be fully quantized too.

\bigskip

\centerline{
{\bf Table V.} DD algebras with
$\mathfrak{g}=\{\mathfrak{s}_3(\mu),\mathfrak{s}_3(0),\mathfrak{n}_3\}$. }

\footnotesize
\begin{center}
\begin{tabular}{|c|c|c|c|c|}
\hline
$(\mathfrak{g},\mathfrak{g}^\ast)$  &  $(\mathfrak{s}_3(\mu),\mathfrak{n}_3)$  &  
$(\mathfrak{s}_3(\mu),\mathfrak{s}_3(1/\mu))$  & 
$(\mathfrak{s}_3(0),\mathfrak{n}_3)$   &  
$(\mathfrak{n}_3,\mathfrak{n}_3)$   
  \\
\hline
DD &  $\mathfrak{s}_3(\mu)\odot R^3$  &   $\mathfrak{so}(1,3)$ 
& 
$\mathfrak{r}_6$  &  $\mathfrak{n}_5\oplus R$\\
\hline
 \cite{gomez} Nr. &  $15$  &   $16$  &   $15'$  & $17$ \\
\hline \hline
 $[x^0,x^1]$ &  $\mu x^1 - x^2$  &  $\mu x^1 - x^2$ &  $-x^2$ &  $0$\\
\hline
 $[x^0,x^2]$ &  $x^1 + \mu x^2$  &  $x^1 + \mu x^2$ & $x^1$ & $0$\\
\hline
 $[x^1,x^2]$ &  $0$ &  $0$  &   $0$ &  $x^0$\\
\hline
 $[X_0,X_1]$ &  $0$  & $ -\lambda X_1/\mu + \lambda X_2$  & $0$ &  $ \omega
X_2$\\
\hline
 $[X_0,X_2]$ &  $0$  & $-\lambda X_1 - \lambda X_2/\mu$  &  $0$ & $0$\\
\hline
 $[X_1,X_2]$ & $\mu\omega X_0$  &  $0$ & $\omega X_0$  &
$0$\\
\hline
 $[x^0,X_0]$ & $0$  &  $0$ & $0$ & $0$\\
\hline
 $[x^0,X_1]$ & $\mu\omega x^2 - \mu X_1 - X_2$  & $-\mu X_1 - X_2$  & $\omega x^2 - X_2$ &
$0$\\
\hline
 $[x^0,X_2]$ & $-\mu\omega x^1 + \mu X_1 - \mu X_2$   &  $X_1 - \mu X_2$  &
$-\omega x^1 + X_1$ & $0$\\
\hline
 $[x^1,X_0]$ & $0$  & $-\lambda x^1/\mu - \lambda x^2$  &  $0$ & $- X_2$\\
\hline
 $[x^1,X_1]$ & $\mu X_0$  & $\lambda x^0/\mu + \mu X_0$  & $0$  &
$0$\\
\hline
 $[x^1,X_2]$ & $-X_0$  & $\lambda x^0 - X_0$   & $-X_0$    & $0$ 
\\
\hline
 $[x^2,X_0]$ &  $0$  & $\lambda x^1 - \lambda x^2/\mu$   &  $0$ & $\omega x^1  +
X_1$
\\
\hline
 $[x^2,X_1]$ & $X_0$  & $-\lambda x^0 + X_0$  & $X_0$  & $-\omega x^0$\\
\hline
 $[x^2,X_2]$ & $\mu X_0$  & $\lambda x^0/\mu + \mu X_0$  & $0$  &$0$ 
\\
\hline
\end{tabular}
\end{center}

\normalsize

\subsection{Case $(\mathfrak{s}_3(\mu),\mathfrak{n}_3)$}

Coproduct:
\[\begin{array}{l}
\Delta(x^0)=1\otimes x^0 + x^0\otimes 1 + \,\mu\omega\,z (x^1\otimes x^2 -
x^2\otimes x^1),\\[0.25cm]  
\Delta(x^1)=1\otimes x^1 + x^1\otimes 1 ,\\[0.25cm]  
\Delta(x^2)=1\otimes x^2 + x^2\otimes 1,\\[0.25cm]  
\Delta(X_0)=1\otimes X_0 + X_0\otimes 1,\\[0.25cm]  
 \Delta(X_1)=e^{-\mu\,z\,X_0}\,\cos(z X_0)\otimes X_1 + X_1\otimes 
e^{\mu\,z\,X_0}\,\cos(z X_0)  \\[0.15cm] 
\hskip2cm
- e^{-\mu\,z\,X_0}\,\sin(z X_0)\otimes X_2 + X_2\otimes 
e^{\mu\,z\,X_0}\,\sin(z X_0),\\[0.25cm]  
 \Delta(X_2)=e^{-\mu\,z\,X_0}\,\cos(z X_0)\otimes X_2 + X_2\otimes 
e^{\mu\,z\,X_0}\,\cos(z X_0)  \\[0.15cm] 
\hskip2cm
+ e^{-\mu\,z\,X_0}\,\sin(z X_0)\otimes X_1 - X_1\otimes 
e^{\mu\,z\,X_0}\,\sin(z X_0).
\end{array}\]

\noindent Commutation rules:

\noindent a) Subalgebras:
$$ \begin{array}{ll}
 [x^0,x^1]=\mu\,x^1 - x^2, \qquad  & [X_0,X_1]=0, \\[0.25cm]
 [x^0,x^2]=x^1 +\mu\, x^2, \qquad  & [X_0,X_2]=0,  \\[0.25cm]
 [x^1,x^2]=0, \qquad  & [X_1,X_2]=\frac{\omega}{2} \frac{\sinh(2\mu\,z X_0)}{z}.     
\end{array}
$$ 

\noindent b) Crossed relations:
$$ \begin{array}{l}
 [x^0,X_0]=0, \\[0.25cm]
 [x^0,X_1]=-\mu\,X_1 - X_2 - \,\mu\,\omega\,x^1\,\sinh(\mu \, z \, X_0) \sin(z\,X_0) +
\mu\omega\,x^2\,\cosh(\mu\,z\,X_0) \cos(z\,X_0), \\[0.25cm]
[x^0,X_2]=X_1 - \mu\,X_2 - \,\mu\,\omega\,x^1\,\cosh(\mu \, z \, X_0) \cos(z\,X_0) -
\mu\omega\,x^2\,\sinh(\mu\,z\,X_0) \sin(z\,X_0), \\[0.25cm]
 [x^1,X_0]= 0, \hskip7.04cm [x^2,X_0]=0,  \cr
 [x^1,X_1]=\frac{\sinh(\mu z X_0)}{z} \cos(z X_0), \hskip4cm
 [x^2,X_1]=\cosh(\mu z X_0) \frac{\sin(z X_0)}{z},  \\[0.25cm]
 [x^1,X_2]=-\cosh(\mu z X_0) \frac{\sin(z X_0)}{z}, \hskip3.55cm
  [x^2,X_2]=\frac{\sinh(\mu z X_0)}{z} \cos(z X_0).     
\end{array}$$ 

\subsection{Case $(\mathfrak{s}_3(\mu),\mathfrak{s}_3(1/\mu))$}

This DD algebra is self-dual for $\mu=1$ and is isomorphic to $\mathfrak{so}(1,3)$ as an
algebra.

\noindent Coproduct:
$$ \begin{array}{l}
 \Delta(x^0)=1\otimes x^0 + x^0\otimes 1 ,\\[0.25cm]
  \Delta(x^1)=e^{-\frac{\lambda}{ \mu}\,z\,x^0}\,\cos(z \lambda x^0)\otimes x^1 +
x^1 \otimes e^{\frac{\lambda}{ \mu}\,z\,x^0}\,\cos(z \lambda x^0)\\[0.15cm]
 \hskip2cm -  e^{-\frac{\lambda}{ \mu}\,z\,x^0}\,\sin(z \lambda x^0)\otimes
x^2 + x^2 \otimes e^{\frac{\lambda}{ \mu}\,z\,x^0}\,\sin(z \lambda x^0) ,\\[0.25cm]  
 \Delta(x^2)=e^{-\frac{\lambda}{ \mu}\,z\,x^0}\,\cos(z \lambda x^0)\otimes x^2 +
x^2 \otimes e^{\frac{\lambda}{ \mu}\,z\,x^0}\,\cos(z \lambda x^0)\\[0.15cm]
\hskip2cm +  e^{-\frac{\lambda}{ \mu}\,z\,x^0}\,\sin(z \lambda x^0)\otimes
x^1 - x^1 \otimes e^{\frac{\lambda}{ \mu}\,z\,x^0}\,\sin(z \lambda x^0) ,\\[0.25cm]  
 \Delta(X_0)=1\otimes X_0 + X_0\otimes 1, \\[0.25cm]
 \Delta(X_1)=e^{-\mu\,z\,X_0}\,\cos(z X_0)\otimes X_1 + X_1\otimes 
e^{\mu\,z\,X_0}\,\cos(z X_0)  \\[0.15cm]
\hskip2cm 
- e^{-\mu\,z\,X_0}\,\sin(z X_0)\otimes X_2 + X_2\otimes 
e^{\mu\,z\,X_0}\,\sin(z X_0),\\[0.25cm]
 \Delta(X_2)=e^{-\mu\,z\,X_0}\,\cos(z X_0)\otimes X_2 + X_2\otimes 
e^{\mu\,z\,X_0}\,\cos(z X_0) \\[0.15cm]
\hskip2cm
+ e^{-\mu\,z\,X_0}\,\sin(z X_0)\otimes X_1 - X_1\otimes 
e^{\mu\,z\,X_0}\,\sin(z X_0).
\end{array}$$ 

\noindent Commutation rules:

\noindent a) Subalgebras:
$$ \begin{array}{ll}
 [x^0,x^1]=\mu\,x^1 - x^2, \qquad  & [X_0,X_1]=-\frac{\lambda}{\mu}\,X_1 + \lambda\,
X_2,  \\[0.25cm]
 [x^0,x^2]=x^1 + \mu x^2, \qquad  & [X_0,X_2]=-{\lambda}\,X_1 -\frac{\lambda}{\mu}\,
X_2,   \\[0.25cm]
 [x^1,x^2]=0, \qquad  & [X_1,X_2]=0.     
\end{array}
$$ 

\noindent b) Crossed relations:
$$ \begin{array}{l}
 [x^0,X_0]=0, \hskip4.96cm   [x^1,X_0]= -\frac{\lambda}{\mu}\,x^1 - \lambda\, x^2, 
\\[0.25cm] 
 [x^0,X_1]=-\mu\,X_1 - X_2, \hskip3cm
[x^1,X_1]=\frac{\sinh(z(\tfrac{\lambda}{\mu} x^0 + \mu\,X_0) )}{z} 
\cos(z ({\lambda} x^0 - \,X_0)),   \\[0.25cm] 
 [x^0,X_2]=X_1 - \mu\,X_2, \hskip3.32cm
  [x^1,X_2]= \frac{\sin(z ({\lambda} x^0 -
\,X_0) )}{z} \cosh(z(\tfrac{\lambda}{\mu} x^0 + \mu\,X_0) ),\\[0.25cm]    
  [x^2,X_0]=\lambda\,x^1 - \frac{\lambda}{\mu}\,x^2,\\[0.25cm] 
[x^2,X_1]= - \frac{\sin(z ({\lambda} x^0 -
\,X_0) )}{z} \cosh(z(\tfrac{\lambda}{\mu} x^0 + \mu\,X_0) ),\\[0.25cm]    
 [x^2,X_2]=\frac{\sinh(z(\tfrac{\lambda}{\mu} x^0 + \mu\,X_0) )}{z} \cos(z ({\lambda} x^0
-\,X_0)).
\end{array}$$

\subsection{Case $(\mathfrak{s}_3(0),\mathfrak{n}_3)$}

Coproduct:
$$ \begin{array}{l}
\Delta(x^0)=1\otimes x^0 + x^0\otimes 1 + \,\omega\,z (x^1\otimes x^2 -
x^2\otimes x^1) ,\\[0.25cm]   
\Delta(x^1)=1\otimes x^1 + x^1\otimes 1,\\[0.25cm]  
 \Delta(x^2)=1\otimes x^2 + x^2\otimes 1,\\[0.25cm]  
 \Delta(X_0)=1\otimes X_0 + X_0\otimes 1,\\[0.25cm] 
 \Delta(X_1)=\cos(z X_0)\otimes X_1 + X_1\otimes  \cos(z X_0)    
-\sin(z X_0)\otimes X_2 + X_2\otimes  \sin(z X_0),\\[0.25cm] 
 \Delta(X_2)= \cos(z X_0)\otimes X_2 + X_2\otimes 
 \cos(z X_0) +  \sin(z X_0)\otimes X_1 - X_1\otimes  \sin(z X_0),
\end{array}$$

\noindent Commutation rules:

\noindent a) Subalgebras:
$$ \begin{array}{ll}
 [x^0,x^1]=  - x^2, \qquad  & [X_0,X_1]=0,  \\[0.25cm] 
 [x^0,x^2]=x^1 , \qquad  & [X_0,X_2]=0,   \\[0.25cm] 
 [x^1,x^2]=0, \qquad  & [X_1,X_2]=\omega   X_0.     
\end{array}$$ 

\noindent b) Crossed relations:
$$ \begin{array}{lll}
 [x^0,X_0]=0, \qquad  & [x^1,X_0]= 0, \qquad & [x^2,X_0]=0,\\[0.25cm] 
 [x^0,X_1]=- X_2  +
 \omega\,x^2\,\cos(z\,X_0), \qquad  & [x^1,X_1]=0, \qquad & [x^2,X_1]=  \frac{\sin(z
X_0)}{z}, \\[0.25cm] 
 [x^0,X_2]=X_1  -\omega\,x^1\, \cos(z\,X_0) , \qquad  & [x^1,X_2]=-  \frac{\sin(z
X_0)}{z},  \qquad &  [x^2,X_2]= 0.
\end{array}
$$ 

\subsection{Case $(\mathfrak{n}_3,\mathfrak{n}_3)$}

This is an essentially self-dual DD algebra, whose quantization is:

\noindent Coproduct:
$$ \begin{array}{l}
 \Delta(x^0)=1\otimes x^0 + x^0\otimes 1 ,\\[0.25cm] 
 \Delta(x^1)=1\otimes x^1 + x^1\otimes 1 ,\\[0.25cm] 
 \Delta(x^2)=1\otimes x^2 + x^2\otimes 1 + 
{z\omega}( x^0\otimes x^1 - x^1\otimes x^0),\\[0.25cm]  
 \Delta(X_0)=1\otimes X_0 + X_0\otimes 1 - z(X_1\otimes
X_2 -  X_2\otimes  X_1),\\[0.25cm]  
\Delta(X_1)=1\otimes X_1 + X_1\otimes 1,\\[0.25cm] 
 \Delta(X_2)=1\otimes X_2 + X_2\otimes 1.
\end{array}$$

\noindent a) Subalgebras:

$$ \begin{array}{ll}
 [x^0,x^1]=0, \qquad  & [X_0,X_1]= \omega X_2,  \\[0.25cm] 
 [x^0,x^2]=0, \qquad  & [X_0,X_2]= 0,  \\[0.25cm] 
 [x^1,x^2]=x^0, \qquad  & [X_1,X_2]=0.     
\end{array}
$$ 

\noindent b) Crossed relations:
$$ \begin{array}{lll}
 [x^0,X_0]=0, \qquad  & [x^1,X_0]= - X_2 ,\qquad  &[x^2,X_0]=X_1 + \omega x^1  \\[0.25cm] 
 [x^0,X_1]=0, \qquad  & [x^1,X_1]=0, \qquad  & [x^2,X_1]=- {\omega}\,x^0\\[0.25cm] 
 [x^0,X_2]=0, \qquad  & [x^1,X_2]=0, \qquad  &  [x^2,X_2]=0.     
\end{array}$$ 

\section{Concluding remarks}

We have presented a comprehensive study of the explicit Hopf algebra quantization of
6D Drinfel'd doubles, thus obtaining a relevant set of new 6D quantum algebras. The only
four cases for which we have not succeeded in the quantization belong to the family of
classical doubles with
$\mathfrak{g}=\mathfrak{r}_3(1)$, and could be obtained by applying specific twists on
the well-known standard quantization of $\mathfrak{so}(1,3)$ and $\mathfrak{sl}_2\odot
R^3$ \cite{so22}. In general, it becomes
apparent that the complexity of the quantum commutation rules that we have obtained is
mainly encoded --by construction-- in the deformation of the crossed relations
$[x^i,X_j]$ within the DD algebra.

This study shows that the quantization of the canonical Lie bialgebra structure of a
classical double can be directly addressed without making use of the cumbersome
construction of the  universal
${\cal T}$-matrix, which is defined on the canonical dual of the quantum universal
enveloping algebra.  Obviously, our ``analytic" approach does not give any
hint concerning the explicit form of the quantum $R$-matrices that should be associated
to each of the quantum doubles here constructed. Certainly, this is an open problem.

Finally we stress that, from the point of view of Drinfel'd doubles, non-simple Lie
algebras should play a relevant role 
for the classification of quantum groups.

\section*{Acknowledgments}

This work was partially supported  by the Ministerio de
Educaci\'on y Ciencia   (Projects BMF2002-02000 and FIS2004-07913),  by the Junta
de Castilla y Le\'on   (Projects VA085/02 and BU04/03), and by INFN-CICyT (Italy-Spain).


\end{document}